\documentclass{amsart}
\usepackage{amsmath,amsfonts,amssymb}

\newtheorem{theorem}{Theorem}[section]

\theoremstyle{definition}

\theoremstyle{remark}

\numberwithin{equation}{section}

\begin{document}

\title[A new proof of the discrete Hardy inequality]{A new proof of the discrete Hardy inequality}

\author[Ivan Gadjev]{Ivan Gadjev}
\address{Department of Mathematics and Informatics,
University of Sofia,
5 James Bourchier Blvd., 
1164 Sofia, Bulgaria}
\email{gadjev@fmi.uni-sofia.bg}

      

\begin{abstract}
We present a short, self-contained elementary proof of the discrete Hardy inequality
$$
\sum_{k=1}^{\infty}\Big(\frac{1}{k}\sum_{j=1}^{k}a_j\Big)^p\leq
\left(\frac{p}{p-1}\right)^{p}\,\sum_{k=1}^{\infty}{a_k^p}
$$
for $p>1$ and every sequence $\{a_k\}_{k=1}^{\infty}$ of non-negative real numbers, along with a direct characterization of the case of equality.
\end{abstract}

\maketitle

\section{Introduction and main result}

The discrete Hardy inequality is one of the cornerstones of classical analysis, with far-reaching applications in harmonic analysis, operator theory, and partial differential equations. Motivated by his investigations into simpler proofs of Hilbert's double-series theorem, G.\,H.~Hardy \cite{H1919, H1920, H1925} proved that for every $p>1$ and every sequence $\{a_k\}_{k=1}^{\infty}$ of non-negative real numbers,
\begin{equation}
\label{eq02}
\sum_{k=1}^{\infty}\Big(\frac{1}{k}\sum_{j=1}^{k}a_j\Big)^p\leq
\left(\frac{p}{p-1}\right)^{p}\,\sum_{k=1}^{\infty}{a_k^p},
\end{equation}
whenever the series on the right-hand side converges.

In his 1920 paper \cite{H1920}, Hardy established \eqref{eq02} with the non-optimal constant $(p^2/(p-1))^p$, noting that Marcel Riesz had independently obtained the same bound. Shortly thereafter, in a 1921 letter to Hardy \cite{Lan2} (later published in \cite{Lan1}), E.~Landau proved \eqref{eq02} with the sharp constant $(p/(p-1))^p$, demonstrating via the test sequence $a_k^\ast = k^{-1/p-\varepsilon}$ as $\varepsilon \to 0^+$ that no smaller constant is admissible. Landau also observed that equality in \eqref{eq02} holds if and only if $a_k = 0$ for all $k \in \mathbb{N}$.

Over the decades, several elementary proofs of \eqref{eq02} have been found, most notably the summation-by-parts arguments of Elliott \cite{Ell1926} and Broadbent \cite{Bro1937} (see also \cite[Theorem~326]{HLP}). The purpose of this note is to provide an alternative, self-contained proof that avoids summation by parts altogether, relying solely on H\"older's inequality with telescoping weights and an elementary application of Jensen's inequality.

Our main result can be stated as follows.

\begin{theorem}\label{thm:main}
Let $p > 1$ and let $\{a_k\}_{k=1}^{\infty}$ be a sequence of non-negative real numbers such that $\sum_{k=1}^{\infty} a_k^p < \infty$. Then \eqref{eq02} holds, with equality if and only if $a_k = 0$ for every $k \in \mathbb{N}$.
\end{theorem}

\section{Proof of Theorem~\ref{thm:main}}

\begin{proof}[Proof of Theorem~\ref{thm:main}]
Let $q = \frac{p}{p-1} > 1$ be the conjugate exponent of $p$, so that $\frac{1}{p} + \frac{1}{q} = 1$. We first establish the finite version of inequality~\eqref{eq02}, namely, for any $n \in \mathbb{N}$ and any non-negative numbers $a_1, \dots, a_n$,
\begin{equation}
\label{eq_fin}
\sum_{k=1}^{n}\Big(\frac{1}{k}\sum_{j=1}^{k}a_j\Big)^p\leq
\left(\frac{p}{p-1}\right)^{p}\,\sum_{k=1}^{n}{a_k^p}.
\end{equation}

By H\"older's inequality, for any two sequences $\mu_i > 0$ and $\eta_i \geq 0$, $i=1, \dots, n$, we have
\[
\sum_{i=1}^{k}\mu_i\eta_i\le\left(\sum_{i=1}^{k}\eta_i^p\right)^{1/p}\left(\sum_{i=1}^{k}\mu_i^q\right)^{1/q},
\]
or equivalently,
\[
\left(\frac{1}{k}\sum_{i=1}^{k}\mu_i\eta_i\right)^p 
\le \frac{1}{k^p}\left(\sum_{i=1}^{k}\eta_i^p\right)\left(\sum_{i=1}^{k}\mu_i^q\right)^{p/q}.
\]
Setting $a_i = \mu_i\eta_i$ (so that $\eta_i = a_i / \mu_i$), summing over $k = 1, \dots, n$, and changing the order of summation yields
\begin{equation}
\label{eq_sum_M}
\sum_{k=1}^{n}\left(\frac{1}{k}\sum_{i=1}^{k}a_i\right)^p
\le \sum_{i=1}^{n} M_i a_i^p \le \left( \max_{1\le i\le n} M_i\right) \sum_{i=1}^{n} a_i^p,
\end{equation}
where 
\[
M_i = \frac{1}{\mu_i^p} M_i^*, \qquad M_i^* = \sum_{k=i}^{n}\frac{1}{k^p}\left(\sum_{j=1}^{k}\mu_j^q\right)^{p/q}.
\]

Choose the telescoping weights $\mu_j = \left(j^{1/q} - (j-1)^{1/q}\right)^{1/q} > 0$ for $j = 1, \dots, n$. Then
\[
\sum_{j=1}^{k}\mu_j^q = \sum_{j=1}^{k} \left(j^{1/q} - (j-1)^{1/q}\right) = k^{1/q}.
\]
Substituting this into $M_i^*$ and using $\frac{p}{q} \cdot \frac{1}{q} - p = -\left(1 + \frac{1}{q}\right)$, we obtain
\[
M_i^* = \sum_{k=i}^{n} \frac{1}{k^{1+1/q}}.
\]
Since the function $t \mapsto t^{-1-1/q}$ is strictly convex on $(0, \infty)$, Jensen's inequality gives
\[
\frac{1}{k^{1+1/q}} < \int_{k-1/2}^{k+1/2} t^{-1-1/q}\,dt = q\left(\left(k - \frac{1}{2}\right)^{-1/q} - \left(k + \frac{1}{2}\right)^{-1/q}\right)
\]
for every $k \ge 1$. Summing from $k = i$ to $n$ yields
\[
M_i^* < q\left(\left(i - \frac{1}{2}\right)^{-1/q} - \left(n + \frac{1}{2}\right)^{-1/q}\right) < q\left(i - \frac{1}{2}\right)^{-1/q}.
\]
Similarly, since $t \mapsto t^{-1/p}$ is strictly convex on $(0, \infty)$, Jensen's inequality gives for every $i \ge 1$:
\[
\mu_i^q = i^{1/q} - (i-1)^{1/q} = \frac{1}{q} \int_{i-1}^{i} t^{-1/p}\,dt > \frac{1}{q}\left(i - \frac{1}{2}\right)^{-1/p},
\]
and hence $\mu_i^p > q^{-p/q}\left(i - \frac{1}{2}\right)^{-1/q}$. Combining the upper bound for $M_i^*$ and the lower bound for $\mu_i^p$, we conclude that for every $1 \le i \le n$,
\[
M_i = \frac{M_i^*}{\mu_i^p} < \frac{q\left(i - \frac{1}{2}\right)^{-1/q}}{q^{-p/q}\left(i - \frac{1}{2}\right)^{-1/q}} = q^{1 + p/q} = q^p = \left(\frac{p}{p-1}\right)^p.
\]

This proves \eqref{eq_fin}. Letting $n \to \infty$ in \eqref{eq_fin} yields \eqref{eq02}. Moreover, if $\{a_k\}_{k=1}^{\infty}$ is not identically zero, then $a_i > 0$ for at least one index $i \ge 1$. Choosing $n \ge i$ sufficiently large so that $\sum_{k=1}^{n} a_k^p > 0$, we have $$\sum_{k=1}^{n}\left(\frac{1}{k}\sum_{j=1}^{k}a_j\right)^p < \left(\frac{p}{p-1}\right)^p \sum_{k=1}^{n} a_k^p$$ strictly by \eqref{eq_sum_M}, and passing to the limit as $n \to \infty$ preserves the strict inequality because the deficit $\left(\frac{p}{p-1}\right)^p a_i^p - M_i a_i^p > 0$ is bounded away from zero independently of $n$. Hence equality in \eqref{eq02} holds if and only if $a_k = 0$ for every $k \in \mathbb{N}$.
\end{proof}

\bibliographystyle{amsplain}

\end{document}